\documentclass[11pt, a4paper,]{article}
\usepackage{ amsfonts,amsmath}
\usepackage[francais, english]{babel}
\newtheorem{prop}{Proposition}[section]
\newtheorem{coro}{Corollary}[section]

\newcommand{\nab}{\nabla}
\newcommand{\anabla}{{\mathcal{A}}^{\nabla}}
\newcommand{\na}{\stackrel{\star}\nabla}
\newcommand{\sn}{\stackrel{\star}{A}_{\xi}}
\newtheorem{defi}{Definition}[section]
\newcommand{\ov}{\overline}
\newtheorem{theo}{Theorem}[section]

\newtheorem{rem}{Remark}[section]
\newtheorem{fact}{Fact}[section]

\title{${\mathcal{A}}^{\nabla}$-tensors on ligktlike hypersurfaces }
\author{C. Atindogbe\footnote{atindogb@iecn.u-nancy.fr,~  Permanent adress: Institut de mathematiques et de sciences Physiques (IMSP), Universit\'e d'Abomey-Calavi(UAC), Benin), 01 BP 613 Porto-Novo, Benin. Email: atincyr@imsp-uac.org} \qquad \qquad L. Berard-Bergery\footnote{berard@iecn.u-nancy.fr}\cr Institut Elie Cartan, Universit\'e Henri Poincar\'e, Nancy I, B.P. 239\cr 54506 Vand\oe uvre-l\`es Nancy Cedex, France}
\date{.}
\begin{document}
\maketitle

\begin{abstract}
\noindent
This paper introduces $\anabla$-tensors on lightlike hypersurfaces $M^{n+1}$ of signature $(0,n)$, $(n\geq 1)$ and investigates on their properties in connection with the null geometry of $M$. In particular, we show that there is an interplay between existence of $\anabla$-tensors of certain type and lightlike warped product structures.
\end{abstract}


\vspace{0.5cm}

\section{Introduction }
\label{intro}
Natural linear conditions generalizing Einstein metric equation are discussed in \cite{Bes} and illustraded by interesting examples. Among such generalizations are $\mathcal{A}$-manifolds (introduced by A. Gray \cite{Gray}), that is Riemannian manifolds $(M,g)$ whose Ricci tensor $r$ satisfies $\nabla r (X,X,X) = 0$ for all $X ~\in TM$, where $\nabla$ is the Levi-Civita connection of the metric $g$. Examples of compact manifolds of this type,other than Einstein or locally products are compact quotients of naturally reductive homogeneous Riemannian manifolds and nilmanifolds covered by the generalized Heisenberg group of A. Kaplan (see \cite{Bes} and references therein). Also, W. Jelonek in \cite{Jel1} gives explicit examples of compact non-homogeneous proper complete $\mathcal{A}$-manifolds, and an example of locally non-homogeneous proper complete one.

A natural generalization of $\mathcal{A}$-manifolds condition is in considering on the Riemannian manifold $(M,g)$ a symmetric $(0,2)$ tensor $\phi$ (or equivalently, since $g$ is non-degenerate, a symmetric tensor $S\in End(TM)$) satisfying an additional condition $\nabla \phi(X,X,X) = 0$. Such tensors are considered and studied in $\cite{Jel1, Jel2, Jel3}$ and called $\mathcal{A}$-tensors (or Killing tensor for $\phi$). In particular, a description of compact Einstein-Weyl manifolds is given in \cite{Jel3} in terms of these tensors. The present paper aims to investigate similar tensors, namely $\anabla$-tensors on lightlike hypersurfaces, in connection with the null geometry of the latter. 

As it is well known, contrary to timelike and spacelike hypersurfaces, the geometry of a lightlike hypersurface $M$ is different and rather difficult since the normal bundle and the tangent bundle have non-zero intersection. At each point $x\in M$, a straight line orthogonal to $M$ lies in $T_{x}M$ and the familly of these straight lines does not determine a normalization of $M$ and consequently an affine connection on $M$. To overcome this difficulty, a theory on the differential geometry of lightlike hypersurfaces developed by Duggal and Bejancu \cite{DB} introduces a non-degenerate screen distribution and construct the corresponding lightlike transversal vecor bundle. This enable to define an induced linear connection (depending on the screen distribution, and hence is not unique in general). On the other hand, it is important to notice that the second fundamental form is independant of the choice of the screen distribution.
 
We brief in section~\ref{prelim} basic informations on normalizations\cite{DB} and pseudo-inversion of degenerate metrics~\cite{ATE}. Our approach in studying $\anabla$-tensors comes from an adaptation of techniques in \cite{Jel1, Jel2} to the case of lightlike hypersurfaces. A known  important result on lightlike hypersurfaces (Theorem~\ref{theo2} below) states that the induced connection is independant of the screen distribution if and only if the lightlike hypersurface is totally geodesic. Equivalently, the induced  connection is torsion-free and metric. In this respect, we introduce in section~\ref{anab}, $\anabla$-tensor (Definition~\ref{ana}) on totally geodesic lightlike hypersurfaces endowed with a specific given screen distribution $S(TM)$ where $\nabla$ is then the unique induced connection on $(M,g)$ in $(\ov{M},\ov{g})$. Thereafter, we show a technical result on its characterisation (Proposition~\ref{caract}). Section~\ref{examp} is concerned  with some explicit constructions (examples) of such tensors. In section~\ref{facts}, we study some geometric properties of these tensors and in section~\ref{3eigenbis} we  establish for a totally geodesic screen distribution, necessary and sufficient condition for eigenspace distributions of $\anabla$-tensors with exactly three eigenspaces to be  integrables(Theorem~\ref{3eigen}). Section~\ref{totumb} is devoted to the special case of totally umbilical screen foliation. In section~\ref{prodstr} we establish a sufficient condition for $\anabla$-tensors to be isotropic. Finally, we show in section~\ref{warped} that there is an interplay between existence of $\anabla$-tensors of certain type and lightlike warped product structure.

\section{Preliminaries on Lightlike hypersurfaces}
\label{prelim}
Let $M$ be a hypersurface  of an $(n+2)-$dimensional pseudo-Riemannian manifold $(\ov{M},\ov{g})$~of index $0 < \nu < n+2 $. In the classical theory of nondegenerate hypersurfaces, the normal bundle has trivial intersection $\{0\}$ with the tangent one and plays an important role in the introduction of main geometric objects. In case of lightlike (degenerate, null) hypersurfaces, the situation is totally different. The normal bundle $TM^{\perp}$ is a rank-one distribution over $M$: $TM^{\perp}\subset TM$ and then coincide with the so called \emph{radical distribution} $RadTM = TM \cap TM^{\perp}$.  Hence,the induced metric tensor field $g$ is degenerate and has rank $n$. The following characterisation is proved in \cite{DB}.   

\begin{prop}  
\label{prop1}
Let $(M,g)$ be a hypersurface of an $(n+2)-$dimensional pseudo-Riemannian manifold $(\ov{M},\ov{g})$. Then the following assertions are equivalent.
\begin{itemize}
\item[(i)]
$M$ is a lightlike hypersurface of $\ov{M}$.
\item[(ii)]
$g$ has constant rank $n$ on $M$.
\item[(iii)]
$TM^{\perp}= \cup_{x\in M}T_xM^\perp$ is a distribution on $M$.
\end{itemize}
\end{prop}
A complementary  bundle of $TM^\perp$ in $TM$ is a rank $n$ nondegenerate distribution over $M$. It is called a \emph{screen distribution} on $M$  and is often denoted by $S(TM)$. A lightlike hypersurface endowed with a specific screen distribution is denoted by the triple $(M,g,S(TM))$. As $TM^\perp$ lies in the tangent bundle, the following result has an important role in studyng the geometry of a lightlike hypersurface.

\begin{prop}(\cite{DB})
\label{theo1}
Let $(M,g,S(TM))$ be a lightlike hypersurface of $(\ov{M},\ov{g})$ with a given screen distribution $S(TM)$. Then there exists a unique rank $1$ vector subbundle $tr(TM)$ of $T\ov{M}|_{M}$, such that for any non-zero section $\xi$ of $TM^\perp$ on a coordinate neighbourhood ${\mathcal{U}}\subset M$, there exists a unique section $N$ of $tr(TM)$ on ${\mathcal{U}}$ satisfyng 
\begin{equation}
\label{eq1}
\ov{g}(N,\xi)= 1
\end{equation}
and
\begin{equation}
\label{eq2}
\ov{g}(N,N)= \ov{g}(N,W)= 0, \quad \quad \forall~ W \in \Gamma(ST|_{\mathcal{U}}).
\end{equation}
\end{prop}
Here and in the sequel we denote by $\Gamma(E)$ the ${\mathcal{F}}(M)-$module of smooth sections of a vector bundle $E$ over $M$, ${\mathcal{F}}(M)$ being the algebra of smooth functions on $M$. Also, by $\perp$ and $\oplus$ we denote the orthogonal and non-orthogonal direct sum of two vector bundles.  By proposition~\ref{theo1} we may write down the following decompositions.

\begin{equation}
\label{eq3}
TM=S(TM) \perp TM^\perp,
\end{equation}

\begin{equation}
\label{eq4}
T\ov{M}|_{M} = TM \oplus tr(TM)
\end{equation}

and

\begin{equation}
\label{eq4bis}
T\ov{M}|_{M}= S(TM) \perp (TM^\perp \oplus tr(TM))
\end{equation}

As it is well known, we have the following:
\begin{defi}
\label{induced}
Let $(M,g,S(TM))$ be a lightlike hypersurface of $(\ov{M},\ov{g})$ with a given screen distribution $S(TM)$. The induced connection, say $\nabla$, on $M$ is defined by
\begin{equation}
\label{eq48}
\nabla_{X}Y = Q(\ov{\nabla}_XY),
\end{equation} 
where $\ov{\nabla}$ denotes  the Levi-civita connection on $(\ov{M},\ov{g})$ and $Q$ is the projection on $TM$ with respect to the decomposition $(\ref{eq4})$.
\end{defi}

\begin{rem}
\label{rem1}
Notice that the induced connection $\nabla$ on $M$ depends on both $g$ and the specific given screen distribution $S(TM)$ on $M$.
\end{rem}
By respective projections $Q$ and $I-Q$, we have Gauss an Weingarten formulae in the form

\begin{equation}
\label{eq5}
\ov{\nabla}_XY= \nabla_XY + h(X,Y)\qquad \forall X,Y ~ \in \Gamma(TM),
\end{equation}

\begin{equation}
\label{eq6}
\ov{\nabla}_XV= -A_VX + \nabla_X^tV \qquad \forall X ~ \in \Gamma(TM),\quad \forall ~V \in \Gamma(tr(TM)).
\end{equation}
Here, $\nabla_XY$ and $A_{V}X$ belong to $\Gamma(TM)$. Hence

$\bullet$  $h$ is a $\Gamma(tr(TM))$-valued symmetric ${\mathcal{F}}(M)$-bilinear form on $\Gamma(TM)$,

$\bullet$ $A_{V}$ is an  ${\mathcal{F}}(M)$-linear operator on $\Gamma(TM)$, and 

$\bullet$ $\nabla^t$ is a linear connection on the lightlike transversal vector bundle $tr(TM)$. 

Let $P$  denote the projection morphism of $\Gamma(TM)$ on $\Gamma(S(TM))$ with respect to the decomposition (\ref{eq3}). We have 

\begin{equation}
\label{eq9}
\nabla_XPY= \na_XPY + h^{*}(X,PY)\qquad \forall X,Y ~ \in \Gamma(TM),
\end{equation}

\begin{equation}
\label{eq10}
\nabla_X U= -\stackrel{\star}{A}_{U}X + \nabla^{*t}_XU \qquad \forall X ~ \in \Gamma(TM),\quad \forall ~U \in \Gamma(TM^\perp).
\end{equation}

Here $\na_XPY$ and $\stackrel{\star}{A}_{U}X$ belong to $\Gamma(S(TM))$,  $\na$ and $\nabla^{*t}$are linear connection on $S(TM)$ and $TM^\perp$, respectively. Hence

$\bullet$  $h^{*}$ is a $\Gamma(TM^\perp)$-valued  ${\mathcal{F}}(M)$-bilinear form on $\Gamma(TM)\times \Gamma(S(TM))$, and 

$\bullet$ $\stackrel{\star}{A}_{U}$ is a $\Gamma(S(TM))$-valued ${\mathcal{F}}(M)$-linear operator on $\Gamma(TM)$.
 
\noindent
They are the second fundamental form  and the shape operator of the screen distribution, respectively. 

Equivalently, consider a normalizing pair $\{\xi, N\}$ as in the proposition~\ref{theo1}. Then, $(\ref{eq5})$ and $(\ref{eq6})$ take the form

\begin{equation}
\label{eq7}
\ov{\nabla}_XY= \nabla_XY + B(X,Y)N \qquad \forall X,Y ~ \in \Gamma(TM|_{\mathcal{U}}),
\end{equation}
and
\begin{equation}
\label{eq8}
\ov{\nabla}_XN= -A_{N}X+ \tau(X)N \qquad \forall X ~ \in 
\Gamma(TM|_{\mathcal{U}}),
\end{equation}
where we put locally on ${\mathcal{U}}$,
\begin{equation}
\label{eq13}
B(X,Y) = \ov{g}(h(X,Y),\xi)
\end{equation}

\begin{equation}
\label{eq14}
\tau(X) = \ov{g}(\nabla^{t}_XN,\xi)
\end{equation}

\noindent
It is important to stress the fact that the local second fundamental form $B$ in $(\ref{eq13})$ does not depend on the choice of the screen distribution.

We also define (locally) on ${\mathcal{U}}$ the following:
\begin{equation}
\label{eq11}
C(X,PY) = \ov{g}(h^{*}(X,PY),N),
\end{equation}

\begin{equation}
\label{eq12}
\varphi(X) = - \ov{g}(\nabla^{\star t}_X\xi ,N).
\end{equation}

Thus, one has for $~X\in~\Gamma(TM)$
\begin{equation}
\label{eq15}
\nabla_XPY= \na_XPY + C(X,PY)\xi
\end{equation}

\begin{equation}
\label{eq16}
\nabla_X\xi = -\sn X + \varphi(X)\xi 
\end{equation}
It is straighforward to verify that for $~X,Y \in~\Gamma(TM)$
\begin{equation}
\label{eq17}
B(X,\xi) = 0 
\end{equation}

\begin{equation}
\label{eq18}
B(X,Y) = g(\sn X,Y)
\end{equation}

\begin{equation}
\label{eq19}
\sn\xi = 0
\end{equation}

The linear connection $\na$ from (\ref{eq9})is a metric connection on $S(TM)$ and we have for all tangent vector fields $X$, $Y$ and $Z$ in $TM$
\begin{equation}
\label{eq20}
\left(\nabla_{X}g\right)(Y,Z)~=~B(X,Y)\eta(Z) + B(X,Z)\eta(Y).
\end{equation}
with
\begin{eqnarray*}
\label{eq49}
\eta(\cdot) = \ov{g}(N,\cdot).
\end{eqnarray*}

The induced connection $\nabla$ is torsion-free, but not necessarily $g$-metric. Also, on the geodesibility of $M$ the following is known.

\begin{theo}(\cite[p.88]{DB})
\label{theo2}
Let $(M,g,S(TM))$ be a lightlike hypersurface of a pseudo-Riemannian manifold $(\ov{M},\ov{g} )$. Then the following assertions are equivalent:
\begin{itemize}
\item[(i)] $M$ is totally geodesic.
\item[(ii)]$h$ (or equivalently $B$) vanishes identically on $M$.
\item[(iii)] $\stackrel{\star}{A}_{U}$ vanishes identically on $M$, for any $U~\in \Gamma(TM^\perp)$
\item[(iv)]The  connection $\nabla$ induced by $\ov{\nabla}$ on $M$ is torsion-free and metric.
\item[(v)] $TM^\perp$ is a parallel distribution with respect to $\nabla$.
\item[(vi)]$TM^\perp$ is a Killing distribution on $M$.
\end{itemize}
\end{theo}
It turns out that if $(M,g)$ is not totally geodesic, there is no connection that is, at the same time, torsion-free and $g$-metric. But there is no unicity of such a connection in case there is any.

\subsection{Pseudo-inversion of degenerate metrics}
A large class of differential operators in differential geometry is intrinsically defined by means of the dual metric $g^{*}$ on the dual bundle $\Gamma(T^{*}M)$~of $1$-forms on $M$. If the metric $g$ is nondegenerate, the tensor field $g^{*}$ is nothing but the inverse of $g$.We brief here construction of some of these operators in case the metric $g$ is degenerate and refer the reader to \cite{ATE} for more details.

Let $(M,g,S(TM))$ be a lightlike hypersurface and $\{\xi,N  \}$ be a pair of  (null-) vectors given by the normalizing theorem \ref{theo1}. Consider on $M$ the one-form defined by
\begin{equation} 
\label{eq20bis}
\eta(\cdot)~=~ \ov{g}(~N~,~\cdot~)
\end{equation}
For all $X\in \Gamma(TM)$,
$$X = PX + \eta(X)\xi $$ 
and $\eta(X)= 0 $ if and only if $X\in \Gamma(S(TM))$. Now, we define $\flat$ by 

\begin{eqnarray*}
\flat: \Gamma(TM)  &  \longrightarrow & \Gamma(T^{*}M)\cr
                   &                  &               \cr 
                X  &  \longmapsto         &  X^{\flat}
\end{eqnarray*}
such that 
\begin{equation} 
\label{eq21}
X^{\flat} = g(~X~, ~\cdot~) + \eta(X)\eta(~\cdot~)
\end{equation}

Clearly, such a $\flat$ is an  isomorphism of $\Gamma(TM)$ onto $\Gamma(T^{*}M)$, and generalize the usual nondegenerate theory. In the latter case,  $\Gamma(S(TM))$ coincide with 
$\Gamma(TM)$, and as a consequence the $1-$form $\eta$ vanishes identically and the projection morphism P becomes the identity map on $\Gamma(TM)$. We let $\sharp$ denote the inverse of the isomorphism $\flat$ given by (\ref{eq21}). For $X\in \Gamma(TM)$ (resp. $\omega \in T^{*}M $), $X^{\flat}$ (resp. $\omega^{\sharp}$) is called the dual $1-$form of $X$ (resp. the dual vector field of $\omega$) with respect to the degenerate metric $g$. It follows (\ref{eq21}) that if $\omega$ is a $1$-form on $M$, we have for $X\in \Gamma(TM)$

\begin{equation}
\label{eq33}
\omega(X)~=~g(\omega^{\sharp},X) ~+~\omega(\xi)\eta(X)
\end{equation}

Now we introduce the so called associate non degenerate metric $\tilde{g}$ to the degenerate metric $g$ as follows. For $X,Y \in \Gamma(TM)$, define $\tilde{g}$ by 

\begin{equation} 
\label{eq22}
\tilde{g}(X,Y)~=~X^{\flat}(Y) 
\end{equation}
Clearly,$ \tilde{g} $ defines a non degenerate metric on $M$ and play an important role in defining the usual differential operators \emph{gradient}, \emph{divergence}, \emph{laplacian} with respect to degenerate metric $g$ on lightlike hypersurfaces. Also, obseve that  $\tilde{g}$ coincides with $g$ if the latter is not degenerate. The $(0,2)$ tensor field $g^{[~\cdot~,~\cdot~]}$, inverse of $\tilde{g}$ is called \emph{the pseudo-inverse of $g$}. Finally, we state the following result (\cite{ATE}).
\begin{prop}
\label{prop2}
Let $(M,g,S(TM))$ be a lightlike hypersurface of a pseudo-Riemannian $(n+2)$-dimensional  manifold $(\ov{M},\ov{g} )$.We have 
\begin{itemize}
\item[(i)] for any smooth function $f:{\mathcal{U}}\subset M \rightarrow \mathbb{R}$,
\begin{equation}
\label{eq23}
grad^{g}f~=~ g^{[\alpha \beta]}f_{\alpha}\partial_\beta 
\end{equation}
where $f_{\alpha}= \frac{\partial f}{\partial x^{\alpha}}$, $\partial_{\beta}= \frac{\partial}{\partial x^{\beta}}$,~~$\alpha,~\beta= 0,\cdots, n $;
\item[(ii)] For any vector field $X$ on ${\mathcal{U}}\subset M$,
\begin{equation}
\label{eq24}
div^{g}X = \sum_{\alpha,\beta}g^{[\alpha,~\beta]}\tilde{g}(\nabla_{\partial_{\alpha}}X,\partial_{\beta})
\end{equation}
\item[(iii)] for smooth function $f:{\mathcal{U}}\subset M \rightarrow \mathbb{R}$
\begin{equation}
\label{eq25}
\Delta^{g}f~=~ \sum_{\alpha \beta}g^{[\alpha,~\beta]}\tilde{g}(\nabla_{\partial_{\alpha}}grad^{g}f,\partial_{\beta})
\end{equation}
\end{itemize}
where $\{\partial_{0}:=\xi, \partial_{1},\cdots,\partial_{n} \}$ is any quasiorthonormal frame field on $M$ adapted to the decomposition (\ref{eq3}).
\end{prop}

In index free notation, (\ref{eq23}) can be written in the form ~$\tilde{g}(\nabla^{g}f,X)~=~df(X)$ which defines the gradient of the scalar function $f$ with respect to the degenerate metric $g$. With nondegenerate $g$, one has $\tilde{g} = g$ so that $(i)-(iii)$ generalize the usual known formulae to the degenerate set up. 

From now on, unless otherwise stated, the ambiant manifold $(\ov{M},\ov{g})$ has a Lorentzian signature so that all lighlike hypersurfaces considered are of signature $(0,n)$. In particular, it follows that any screen distribution is Riemannian. 
As it is well known (theorem~\ref{theo2}), only totally geodesic lightlike hypersurfaces do have their induced connection metric and torsion-free. In the next section and the remainder of the text, only such lightlike hypersurfaces will be in consideration. We also assume that the null vector field $\xi$ is globally defined on $M$. Respective metrics will be denoted $\langle\cdot,\cdot\rangle$  if no ambiguity occurs.

\section{$\anabla$-tensors}
\label{anab}
\begin{defi}
\label{ana}
Let $(M,g,S(TM))$ be a totally geodesic lightlike hypersurface of $(\ov{M},\ov{g})$, $\nab$  the induced (Levi-civita) connection on $M$. By $\anabla$-tensor on $(M,g,S(TM))$, we mean a screen preserving element $S\in End(TM)$ for which 
\begin{enumerate}
\item[(a)] $\langle SX,Y\rangle = \langle X,SY\rangle$ ~for all $X$, $Y$ in $TM$,
\item[(b)] $X^{\flat}(\nabla S(X,X)) = 0$ ~~\qquad for all $X$ in $TM$,
\end{enumerate}

hold, where $\flat$ denote the duality isomorphism between $TM$ and $TM^{\star}$ with respect to the degenerate metric tensor $g$ and the screen distribution $S(TM)$.
\end{defi}

It should be noticed that screen preserving means $P$ and $S$ commute. One also write $S\in \anabla$ if $S$ is an $\anabla$-tensor. An $\anabla$-tensor is called isotropic if it is $Rad(TM)$-valued, otherwise, it is called a proper $\anabla$-tensor.

Killing tensors on $M$ are symmetric $(0,2)$-tensors, say $\phi$ such that 
\begin{equation}
\label{eq26}
\phi(X,Y)=\langle SX, Y\rangle\quad \forall~ X, Y \in TM
\end{equation}
for some $\anabla$-tensor $S$. Observe that $\phi$ is a degenerate $(0,2)$-tensor since at each $u\in M$ its nullity space $\eta_{\phi|_{u}} \supset RadTM|_{u}$, i.e
\begin{equation}
\label{eq27}
\phi(X,\xi)= \phi(\xi,X)=0 \quad \forall~ X \in TM, ~\forall \xi \in  RadTM.
\end{equation}
It also satisfies
\begin{equation}
\label{eq28}
\nabla \phi (X,X,X)=\langle\nabla S(X,X),X\rangle \quad \forall ~ X \in TM.
\end{equation}
Since $\nabla$ is a metric connection, we have 
\begin{equation}
\label{eq29}
\nabla_{X}\xi = \varphi(X)\xi \quad \forall ~ X \in TM,
\end{equation}
for some global $1-$form $\varphi$ on $M$.
The following proposition is the equivalent to the Riemannian case \cite{Jel1}.

\begin{prop}
\label{caract}
Let $(M,g,S(TM))$ be a totally geodesic lightlike hypersurface, $S$ a symmetric $(1,1)$tensor and $\phi(X,Y)=\langle SX,Y\rangle$ for all $X$ in $TM$. The following assertions are equivalent.
\begin{enumerate}
\item[(a)]$S\in \anabla$.
\item[(b)]For every geodesic $\gamma$ on $(M,g)$, the real valued function $t \longmapsto \phi(\gamma'(t), \gamma'(t))$ is constant on $dom\gamma$ and, if $\gamma$ is a null geodesic, $S(\gamma'(t))$ is parallel along $\gamma$.
\item[(c)]$\Sigma_{cyclic}\nabla_{X}\phi(Y,Z)=-\Sigma_{cyclic}\eta(X)\eta(\nabla S(Y,Z))$.
\end{enumerate}
\end{prop}

\noindent
{\bf{Proof.~}}The equivalence $(a)$ and $(c)$ is immediate using definition of $\flat$ and bipolarization of relation $(b)$ in definition~\ref{ana}.Let us show the equivalence $(a)$ and $(b)$.Assume $(a)$ and consider $\gamma$ a geodesic on $M$. We have 
\begin{eqnarray*}
\frac{d}{dt}\phi(\gamma'(t), \gamma'(t))&=& \nabla_{\gamma'(t)}\phi(\gamma'(t), \gamma'(t))
\end{eqnarray*}
We distinguish two cases: $\gamma$ is a null geodesic or not.

If $\gamma$ is a non null geodesic, from (\ref{eq28}) we have 
\begin{eqnarray*}
\frac{d}{dt}\phi(\gamma'(t), \gamma'(t))&=& \nabla_{\gamma'(t)}\phi(\gamma'(t), \gamma'(t))\cr 
&= & \gamma'(t)^{\flat}(\nabla S(\gamma'(t), \gamma'(t)) )=0,
\end{eqnarray*}
i.e $\phi$ is constant on $dom \gamma$.

If $\gamma$ is a null geodesic, it follows definition of $\phi$ that it vanishes identically on $dom \gamma$. In addition, $\gamma'(t)$ is proportionnal to $\xi$ for all $t$ in $dom \gamma$. Thus, there exists a nowhere vanishing function $t\rightarrow \lambda_{0}(t)$ on $dom \gamma$ such that 
\begin{equation}
\label{eq30}
\nabla S(\gamma'(t), \gamma'(t))= (\lambda_{0}(t))^{2}\nabla S(\xi,\xi) \in RadTM|_{\gamma}.
\end{equation}
Using $(b)$ in definition~\ref{ana}, we have $\eta(\nabla S(\gamma'(t), \gamma'(t)))=0 \quad \forall t \in dom \gamma$. This together with (\ref{eq30}) lead to $\nabla S(\gamma'(t), \gamma'(t))=0 \quad \forall t \in dom \gamma$. Finally, since $\gamma$ is a geodesic, we have $\nabla_{\gamma'(t)}S(\gamma'(t))=0$, and $(b)$ is proved.

Conversely, assume $(b)$ holds and let $X \in T_{x_{0}}M$, $x_{0}\in M$. Consider $\gamma$ a geodesic satisfyng initial conditions $\gamma(0)=x_{0}$ and $\gamma'(0)=X$. One has
\begin{eqnarray*}
X^{\flat}(\nabla S(X,X))= \nabla_{\gamma'(t)}\phi(\gamma'(t), \gamma'(t))|_{t=0}+ \eta(\gamma'(t))\eta(\nabla S(\gamma'(t), \gamma'(t)))|_{t=0}=0
\end{eqnarray*}
i.e $(a)$ is proved and the proof is complete.$\square$

\begin{rem}
\label{rem1bis}
\end{rem}
\begin{enumerate}
\item[(a)]Observe that for $X,~Y$and $Z$ in $S(TM)$, relation $(c)$ in proposition~\ref{caract} reduces to 

\begin{equation}
\label{eq31}
\nabla_{X}\phi(Y,Z)+\nabla_{Y}\phi(Z,X)+\nabla_{Z}\phi(X,Y)=0.
\end{equation}
\item[(b)] Since $M$ has signature $(0,n)$, $n= dim M - 1$, the $\anabla$-tensor $S$ induces by restriction on the nondegenerate (Riemannian) screen distribution $S(TM)$, a $\mathcal{A}$-tensor $S'$ with respect to the (unique) Levi-Civita connection $\na$ induced by $\nabla$ on $S(TM)$. Indeed, $S' \in End(S(TM))$ by screen preserving of $S$ and it is known \cite{Jel1} that in this case, (\ref{eq31}) is equivalent to being $\mathcal{A}$-tensor for $S'$. So, the $\anabla$-tensor $S$ splits as 
\begin{equation}
\label{eq32}
S=S'\circ P + \eta(\cdot)S~\xi
\end{equation}
One can show that, if $\sigma$ is a Riemannian $\mathcal{A}$-tensor on $S(TM)$ and if in addition the screen distribution is totally geodesic in $M$, then, for $\lambda_{0}\in C^{\infty}(M)$ , the $(1,1)$-tensor defined on $M$ by 
\begin{equation}
\label{33}
S=\sigma \circ P + \lambda_{0}\eta(\cdot)\xi
\end{equation}
is an $\anabla$-tensor on $M$, provided $\xi \cdot \lambda_{0}=0$.
\end{enumerate}

\section{Constructions. Examples}
\label{examp}
\begin{enumerate}
\item[(a)] Let $M=\mathbb{L}\times M_{1} \times_{f}M_{2}$ be a totally geodesic lightlike triple warped product hypersurface, with $f$ a smooth positive function on $M_{1}$, $\mathbb{L}$ a (one dimensional) global null curve, $(M_{i},g_{i})$ Riemannian manifolds ($i=1,2$). Since $M$ is totally geodesic, it is possible to use a normalization for which the $1$-form $\tau$ (or equivalently $\varphi$) vanishes identically. Let $\nabla^{i}$ $(i=1,2)$ denote the Levi-Civita connection on $(M_{i},g_{i})$. We have 
$$g= g_{1} + (f\pi_{1}(x))^{2}g_{2}$$ 
and the induced connection $\nabla$ on $M$ is given for $X$, $Y$ tangent to $M'=M_{1} \times_{f}M_{2}$ by  
\begin{eqnarray}
\label{34}
\nabla_{X}Y &=& \nabla^{1}_{X_{1}}Y_{1} +\nabla^{2}_{X_{2}}Y_{2}+\left[X_{1}(\psi)Y_{2}+Y_{1}(\psi)X_{2}-g(X_{2},Y_{2})grad\psi\right]\cr & & + C(X,Y)\xi
\end{eqnarray}
where $\pi_{1}$ denotes the projection on the factor $M_{1}$ of $M$,  $X=(X_{1},0)+(0,X_{2})= (X_{1},X_{2})$, $Y=(Y_{1},0)+(0,Y_{2})= (Y_{1},Y_{2})$ on $M_{1}\times M_{2}$, $\nabla^{i}_{X_{i}}Y_{i}|_{p}\in T_{p}M_{i}$ with the vector $(\nabla^{1}_{X_{1}}Y_{1}|_{p},0_{q})\in T_{(p,q)}M_{1}\times M_{2}$ etc., $\psi = \ln f $ and $grad \psi $ its gradient with respect to $g$, and $C$ the second fundamental form of the screen distribution $S(TM) = TM_{1}\oplus TM_{2}$. Note that for $X\in \Gamma(TM)$, due to $[\xi,X]=0$, we have 
\begin{equation}
\label{eq35}
\nabla_{\xi}X = \nabla_{X}\xi = -\tau(X)\xi = 0
\end{equation}
Now, assume that $S(TM)$ is totally geodesic in $M$ (and hence in the ambiant space $\ov{M}\supset M$) and define a $(1,1)$ tensor on $M$ by 
\begin{equation}
\label{eq36}
\left\{
\begin{array}{lcl}
S(\xi)&=&\mu \xi, \quad \mu \in \mathbb{R}\cr
& & \cr
S(X)&=&0, \quad X\in D_{1}=TM_{1}\cr 
& & \cr
S(X)&=&\lambda X, \quad \lambda = Cf^{2}, ~C\in  \mathbb{R}^\star
\end{array}
\right.
\end{equation} 
$S$ is a well defined $(1,1)$ tensor on $M$ that preserves the screen distribution and is obviously symmetric. Let $X=\eta(X)\xi + X_{1} +X_{2} \in TM$. Our aim is to show that $X^{\flat}(\nabla S(X,X))=0$. We have $SX = \mu \eta(X)\xi + \lambda X_{2}$ and direct computation gives
\begin{equation}
\label{37}
\nabla_{X}(SX)=-\lambda g(X_{2},Y_{2})grad\psi + 3Cf^{2}X_{1}(\psi)X_{2} + \lambda\nabla^{2}_{X_{2}}X_{2}
\end{equation}
Also, 
\begin{equation}
\label{38}
S(\nabla_{X}X)= \lambda \left[2X_{1}(\psi)X_{2}+  \nabla^{2}_{X_{2}}X_{2} \right]
\end{equation}
Then, 

\begin{equation}
\label{eq39}
\nabla_{X}(SX)-S(\nabla_{X}X)= \lambda \left[-g(X_{2},X_{2})grad \psi + X_{1}(\psi)X_{2}\right]
\end{equation}
Therefore
\begin{eqnarray*}
X^{\flat}(\nabla S(X,X))&=&X^{\flat}(\nabla_{X}(SX)-S(\nabla_{X}X))\cr
& &\cr
& =&\eta(X)\xi^{\flat}[\nabla_{X}(SX)-S(\nabla_{X}X)]\cr 
& &  +(X_{1}+X_{2})^{\flat}[\nabla_{X}(SX)-S(\nabla_{X}X)]\cr 
& &\cr
&= &\eta(X)\eta[\nabla_{X}(SX)-S(\nabla_{X}X)]\cr & &+(X_{1}+X_{2})^{\flat}[\nabla_{X}(SX)-S(\nabla_{X}X)]\cr 
\end{eqnarray*}
Since by (\ref{eq39}), $\nabla_{X}(SX)-S(\nabla_{X}X)$ is $STM$-valued, we have $$\eta[\nabla_{X}(SX)-S(\nabla_{X}X)=0].$$ But the second term is 
\begin{eqnarray*}
(X_{1}+X_{2})^{\flat}[\nabla_{X}(SX)-S(\nabla_{X}X)]& =& \lambda (X_{1}+X_{2})^{\flat}[-g(X_{2},X_{2})grad \psi \cr & & + X_{1}(\psi)X_{2}]\cr 
& =& \lambda [-\langle X_{1},grad \psi\rangle\langle X_{2},X_{2}\rangle \cr & &+ X_{1}(\psi)\langle X_{2},X_{2}\rangle]=0.
\end{eqnarray*}
Thus, $X^{\flat}(\nabla S(X,X))= 0$ and $S$ defines an $\anabla$-tensor on $(M,g,S(TM))$.

\item[(b)]{\bf{Killing horizons}}.
Let $(M,g)$ be a lightlike hypersurface of a pseudo-Riemannian manifold $(\ov{M},\ov{g})$ and $\ov{G}$ a continuous $k$-parameters group of isommetry acting on $(\ov{M},\ov{g})$. By local isommetry horizon (LIH in short) with respect to $\ov{G}$ is meant a lightlike hypersurface that is invariant under $\ov{G}$ and for which each null geodesic is a trajectory of the group. In case  $\ov{G}$ is $1-$parameter, the LIH is said to be a \emph{killing horizon}. It turns out that a Killing horizon is a lightlike hypersurface whose null tangent vector can be normalized to coincide with a killing vector field \cite{Carter}. Taking into account theorem~\ref{theo2}, killing horizons are totally geodesic in $(\ov{M},\ov{g})$. By \emph{global hypersurface} in a Killing horizon $M$  we mean a topological hypersurface which is crossed exactly once by any null geodesic trajectory of $M$. A killing horizon admitting such a hypersurface will be called \emph{a globally killing horizon}. On the latter, it is possible to construct a special screen distribution as follows. Let $(\varphi_{t})_{t\in I\subset \mathbb{R}}$ be the $1$-parameter group with respect to which $M$ is a killing horizon, and $H$ a global hypersurface in $M$. By definition of $H$ it follows that for each $p\in M$, there exists a unique $(t,q) \in I\times H$ such that $p=\varphi_{t}(q)$. We set $S(T_{p}M)= \varphi_{t\star q}(T_{q}H) $. Clearly, such a $S(TM)$defines an integrable screen distribution on $M$, we denote $S(TM,\varphi_{t},H)$. Recall that throughout the text, the ambiant manifold  $(\ov{M},\ov{g})$ has Lorentzian signature so that global hypersurfaces are Riemannian. Also, the normalized null tangent vector on the killing horizon will be denoted $\xi$. Consider now a globally killing horizon for which local geodesic symmetries preserve a global hypersurface, say $H$, and volume of its regions. The Ricci endomorphism (or the Ricci tensor) of such a $H$ is an $\mathcal{A}$-tensor \cite{Gray}, say $\sigma$. Now define on     $(M,g,S(TM,\varphi_{t},H))$ a $(1,1)$-tensor by 
\begin{equation}
\label{eq40}
SX = \mu\eta(X)\xi + \sigma(PX), \qquad \mu \in \mathbb{R},
\end{equation}
where $P$ denote the projection morphism of the bundle $TM$ on the screen distribution  $S(TM,\varphi_{t},H)$ with respect to the decomposition (\ref{eq3}). Clearly, such a $S$ is $g$-symmetric, preserves $S(TM,\varphi_{t},H)$. Also, observe that since local geodesic symmetries preserve $H$, the screen distribution $S(TM,\varphi_{t},H)$ is totally geodesic in $M$. Finally, using $(b)$ in remark~\ref{rem1bis}, it follows that $S$ is a $\anabla$-tensor on $M$.
\end{enumerate}

\section{Some Facts}
\label{facts}

\begin{fact}
\label{diag}
Any $\anabla$-tensor on $(M,g,S(TM))$ is diagonanalisable.
\end{fact}

\noindent
{\bf{Proof.~}}First, observe that the global null vector $\xi$ is an eigenvector field of $S$. Since $\langle S\xi,X\rangle = \langle \xi,SX\rangle=0$ for all $X$ in $TM$, it follows $S\xi \in RadTM $ and there exists a smooth function $\lambda_{0}$ such that $S\xi = \lambda_{0}\xi$. Since $S(TM)$ is Riemannian, we know that $S'$ is diagonalizable and the same is for $S$ using (\ref{eq32}).$\square$

Now, define the integer-valued function $$x\rightarrow E_{S}(x) = Card\{\mbox{distinct eigenvalues of} S_{x}   \}$$ and set 
$$M_{S}=\{x\in M: E_{S} \mbox{~is constant in a neighbourhood of}~ x \}$$
The set $M_{S}$ is open and dense in $M$. On each component U of $M_{S}$, the dimension, say $p_{\alpha}$,  of the eigenspace $D_{i}=Ker(S-\lambda_{\alpha}I)$ associated to the eigenfunction $\lambda_{\alpha}$ is constant. From now on, we assume all manifolds connected unless otherwise stated and $M=M_{S}$. Also, note that 
\begin{equation}
\label{eq41} 
TM = \sum_{\alpha}^{k}D_{\alpha}
\end{equation}
with $D_{0}= RadTM= span\{\xi  \}$. We use the following range of indices: $0\le \alpha \le k$ and $1\le i \le k$. We have the following technical result.

\begin{fact}
\label{techn}
Let $S$ denote an $\anabla$-tensor on $(M,g,S(TM))$, $\lambda_{0},\lambda_{1},\cdots, \lambda_{k}$ in $C^{\infty}(M)$ be eigenfunctions of $S$. Then,
\begin{eqnarray}
\label{eq42}
\forall~X\in D_{i}, \nabla S(X,X)&=&-\frac{1}{2}\langle X,X\rangle \nabla^{g}\lambda_{i}  \cr 
& & +\left[\frac{1}{2}\langle X,X\rangle \eta(\nabla^{g}\lambda_{i})+(\lambda_{i}-\lambda_{0})C(X,X)  \right]\xi
\end{eqnarray}
and 
\begin{eqnarray}
\label{eq43}
D_{\alpha}&\subset & Ker d\lambda_{\alpha}\qquad 0\le \alpha \le k.
\end{eqnarray}
If $i\ne j$, $X\in \Gamma(D_{i})$ and $Y\in \Gamma(D_{j})$ then 
\begin{eqnarray}
\label{eq44}
\langle\nabla_{X}X,Y\rangle&=&\frac{1}{2} \frac{Y\cdot \lambda_{i}}{\lambda_{j}-\lambda_{i}}\langle X,X\rangle.
\end{eqnarray}

If $X\in \Gamma(D_{0})$ or $Y\in \Gamma(D_{0})$
\begin{eqnarray}
\label{eq45}
\langle \nabla_{X}X,Y\rangle&=&0.
\end{eqnarray}
\end{fact}

\noindent
{\bf{Proof.~}}
For $X\in \Gamma(D_{i})$ and $Y\in \Gamma(TM)$we have
\begin{equation}
\label{eq46}
\nabla S(Y,X)=(Y\cdot \lambda_{i})X + (\lambda_{i}I - S)\nabla_{Y}X.
\end{equation}
Then,
\begin{eqnarray*}
\langle \nabla S(Y,X),X\rangle&=&(Y\cdot \lambda_{i})\langle X,X\rangle + \langle(\lambda_{i}I - S)\nabla_{Y}X,X\rangle \cr
& =&(Y\cdot \lambda_{i})\langle X,X\rangle + \langle \nabla_{Y}X,\lambda_{i}X-\lambda_{i}X\rangle\cr
&=&(Y\cdot \lambda_{i})\langle X,X\rangle
\end{eqnarray*}
that is 
\begin{equation}
\label{eq47}
\langle \nabla S(Y,X),X\rangle = (Y\cdot \lambda_{i})\langle X,X\rangle.
\end{equation}
Therefore, taking $Y=X$ leads to
$$0=\langle\nabla S(X,X),X\rangle = (X\cdot \lambda_{i})\langle X,X\rangle \qquad 1\le i \le k.$$
Since $X\in \Gamma(D_{i})\subset \Gamma(STM)$ (Riemannian), we have $\langle X,X\rangle \ne 0$ and $X\cdot \lambda_{i} = 0 \qquad 1\le i \le k.$ that is $D_{i}\subset Ker d\lambda_{i}  \qquad 1\le i \le k.$
Also, integrale curves of $\xi$ are null geodesics. Then $\nabla S(\xi,\xi)=0 = (\xi \cdot \lambda_{0})\xi$ and  $(\xi \cdot \lambda_{0})= 0$. Thus,  $D_{0}\subset Ker d\lambda_{0}$ and (\ref{eq43}) is proved. From (\ref{eq46}) and (\ref{eq43}) it follows that 
\begin{equation}
\label{eq48bis}
\nabla S(X,X) = (\lambda_{i}I - S)\nabla_X{Y}X.
\end{equation}
Observe that for $X$,$Y$ and  $Z$ in $\Gamma(STM)$, (\ref{eq31}) is equivalent to $$\langle\nabla S(X,Y),Z\rangle + \langle \nabla S(Y,Z),X\rangle + \langle\nabla S(Z,X),Y\rangle =0.$$ Also, $\langle\nabla S(X,Y),X\rangle = \langle\nabla S(X,X),Y\rangle$. Hence $2\langle\nabla S(X,X),Y\rangle + \langle\nabla S(Y,X),X\rangle = 0 $. Taking into account (\ref{eq47}) yields $2\langle\nabla S(X,X),Y\rangle + (Y\cdot \lambda_{i})\langle X,X\rangle = 0$, i.e
\begin{equation}
\label{eq49bis}
\langle2\nabla S(X,X)+\langle X,X\rangle\nabla^{g}\lambda_{i},Y\rangle= 0 \qquad \forall Y\in S(TM).
\end{equation}
Then, since (\ref{eq49bis}) holds trivially for $Y \in RadTM$, 
\begin{equation}
\label{eq50}
\langle2\nabla S(X,X)+\langle X,X\rangle\nabla^{g}\lambda_{i},Y\rangle= 0 \qquad \forall Y\in \Gamma(TM).
\end{equation}
Thus,
$$2\nabla S(X,X)+\langle X,X\rangle\nabla^{g}\lambda_{i} \in RadTM = Span\{\xi\}$$
It follows that 
\begin{equation}
\label{eq51}
\nabla S(X,X)= -\frac{1}{2}\langle X,X\rangle\nabla^{g}\lambda_{i}+ q(X)\xi \qquad \forall Y\in \Gamma(TM).
\end{equation}
where $q(X)$ is a quadratic function in $X$. From (\ref{eq51}), we have 
\begin{equation}
\label{eq52}
\eta(\nabla S(X,X))=-\frac{1}{2}\langle X,X\rangle\eta(\nabla^{g}\lambda_{i})+q(X)
\end{equation}
Now, using (\ref{eq32}), we derive for $X\in \Gamma(D_{i})$,
\begin{equation}
\label{eq53}
\nabla S(X,X))= \na S'(X,X) + (\lambda_{i}-\lambda_{0})C(X,X)\xi
\end{equation}
and 
\begin{equation}
\label{eq54}
\eta\left(\nabla S(X,X)\right)=(\lambda_{i}-\lambda_{0})C(X,X)
\end{equation}
Thus, combining  (\ref{eq52}) and (\ref{eq54}) lead to

\begin{equation}
\label{eq55}
q(X)=\frac{1}{2}\langle X,X\rangle\eta(\nabla^{g}\lambda_{i})+(\lambda_{i}-\lambda_{0})C(X,X).
\end{equation}
Substitute in (\ref{eq51}) to get the announced relation in (\ref{eq42}).

For $X\in \Gamma(D_{i})$, $Y\in \Gamma(D_{j})$ with $i\ne j$,
\begin{eqnarray*}
\langle\nabla S(X,X),Y\rangle&=&\langle(\lambda_{i}I - S)\nabla_{X}X,Y\rangle\cr
                 &=&\langle \nabla_{X}X,(\lambda_{i}-\lambda_{j})Y
		 \rangle.
\end{eqnarray*}
Thus, by (\ref{eq42}),
$$-\frac{1}{2}\langle X,X\rangle\langle\nabla^{g}\lambda_{i},Y\rangle = (\lambda_{i}-\lambda_{j})\langle\nabla_{X}X,Y\rangle$$
and 

\begin{eqnarray*}
\langle\nabla_{X}X,Y\rangle &=& \frac{1}{2}\langle X,X\rangle \frac{Y\cdot \lambda_{i}}{\lambda_{j}-\lambda_{i}}\langle X,X\rangle.
\end{eqnarray*}
Finally, it is clear that if $X\in \Gamma(D_{0})$ ~or $Y\in \Gamma(D_{0})$, one has $\langle \nabla_{X}X,Y\rangle =0$, and the proof is complete.$\square$

\begin{coro}
\label{coro1}
The following assertions are equivalent.
\begin{enumerate}
\item[(a)]$\forall ~ X\in \Gamma(D_{i}), \nabla_{X}X \in \Gamma(D_{i})$.
\item[(b)]$\forall ~ X,Y \in \Gamma(D_{i}), \nabla_{X}Y + \nabla_{Y}X \in \Gamma(D_{i})$.
\item[(c)]$\forall ~ X\in \Gamma(D_{i}), \nabla S(X,X)=0 $.
\item[(d)]$\forall ~ X,Y \in \Gamma(D_{i}), \nabla S(X,Y)+\nabla S(Y,X)=0 $.
\item[(e)]$\nabla^{g}\lambda_{i}$ is $D_{0}$-valued vector field and $\forall~X\in \Gamma(D_{i}),~ C(X,X)=0, \quad 1\le i \le k $.
\end{enumerate}
\end{coro}

\noindent
{\bf{Proof.~}}
The equivalences $(a)\iff (b)$ and  $(c)\iff (d)$  are obvious as polarizations. Let us show $(a)\iff (c)$. We have 
$$\nabla_{X}X \in~ \Gamma(D_{i})~~\stackrel{(\ref{eq48bis})}{\Longrightarrow}~~ \nabla S(X,X)=0.$$
Conversely, if for all $X$ in $\Gamma(D_{i})$, $\nabla S(X,X)=0$, then by (\ref{eq48bis}), $(\lambda_{i}I-S)\nabla_{X}X =0$, i.e $\nabla_{X}X \in \Gamma(D_{i})$, thus $(a)\iff (c)$. Finally, using (\ref{eq42}) we obtain 
$$\nabla S(X,X) = 0 \iff -\frac{1}{2}\langle X,X\rangle P\nabla^{g}\lambda_{i}   +(\lambda_{i}-\lambda_{0})C(X,X)\xi =0 $$
which is equivalent to $P\nabla^{g}\lambda_{i} =0$ and $C(X,X)=0$, i.e $(e)$. This cpmpletes the proof.$\square$ 

Note that $D_{0}$ is of rank one, then is integrable. Also, for $X$, $Y$ in $\Gamma(D_{i})$, we have 
$$\nabla S(X,Y)-\nabla S(Y,X)= (\lambda_{i}I-S)([X,Y])$$ so that $D_{i}$ is  integrable if and only if $\forall ~X,~Y$ in ~$\Gamma(D_{i}),\nabla S(X,Y)-\nabla S(Y,X)$. Moreover, we obtain the following.

\begin{fact}
\label{integr}
If $\nabla^{g}\lambda_{i}$ is $D_{0}$-valued and for all $X\in \Gamma(D_{i})$, $C(X,X)=0$, then the following assertions are equivalent on $M$.
\begin{enumerate}
\item[(a)]$D_{i}$ is integrable.
\item[(b)]For all $X,~Y$ in $\Gamma(D_{i})$, $\nabla S(X,Y)=0$.
\item[(c)]$D_{i}$ is autoparallel.
\end{enumerate}
\end{fact}

\noindent
{\bf{Proof.~}}
For the first equivalence, we shall prove $(a)\Longrightarrow (b) $ and observe that $(b)\Longrightarrow (a)$ is obvious. Assume that $(a)$ holds. From corollary~\ref{coro1}(valid since $(e)$ holds by hypothesis),  $\nabla S(X,Y)+\nabla S(Y,X)=0$ and integrability implies $\nabla S(X,Y)=\nabla S(Y,X)$. Thus,  $\nabla S(X,Y)=0$ and $(a)\Longrightarrow (b) $. Finally, from $\nabla_{X}Y+\nabla_{Y}X \in \Gamma(D_{i})$ and $\nabla_{X}Y-\nabla_{Y}X = [X,Y] \in \Gamma(D_{i})$ we obtain the equivalence $(a)\iff (c)$ .$\square$

\section{$\anabla$-tensors with exactly three eigenspaces}
\label{3eigenbis}
We consider and investigate on some geometric properties  of $\anabla$-tensors with exactly three eigenspaces $D_{0}=Ker(\lambda_{0}I-S)$, $D_{\alpha}=Ker(\alpha I-S)$ and $D_{\beta}=Ker(\beta I-S)$ with $S(TM)=D_{\alpha} \oplus D_{\beta}$. In Riemannian setting, a classical theorem due to Jelonek \cite{Jel1} states that, for a $\mathcal{A}$-tensor with exactly two eigenvalues $\lambda$, $\mu$ and a constant trace, the eigenvalues are necessarily constant, and the eigenspace distributions are both integrable if and only if the $\mathcal{A}$-tensor is parallel. The following is a lightlike version of this result with three eigenvalues.
\begin{theo}
\label{3eigen}
Let $S$ be an $\anabla$-tensor on  $(N,g,S(TM))$ with exactly three eigenfunctions $\lambda_{0} = cte$, $\alpha$, $\beta$ and a constant trace. Then $\nabla^{g}\alpha$  and $\nabla^{g}\beta$ are $D_{0}=RadTM$-valued. In addition, If $S(TM)$ is totally geodesic then the distributions $D_{\alpha}$ and $D_{\beta}$ are both integrable if and only if $\nabla S$ vanishes on $S(TM)\times S(TM)$.
\end{theo}

\noindent
{\bf{Proof.~}}
Since $S$ is smooth, $x\rightarrow p(x)=\mbox{dim} D_{\alpha}(x)$ and $x\rightarrow q(x)=\mbox{dim} D_{\beta}(x)$ are discrete differentiable functions on $M_{S}=M$, so they are constant functions we denote by $p$ and $q$ repectively. From $$\lambda_{0}+p\alpha +q\beta = trS = cte$$
we derive
\begin{equation}
\label{eq56}
p\nabla^{g}\alpha +q\nabla^{g}\beta = 0, \qquad (\nabla^{g}\lambda_{0}=0). 
\end{equation}
Observe that $\langle\nabla^{g}\alpha,\nabla^{g}\beta\rangle = 0$. Then from (\ref{eq56}) we obtain $p\langle\nabla^{g}\alpha,\nabla^{g}\alpha\rangle = 0$ and $q\langle\nabla^{g}\beta,\nabla^{g}\beta\rangle = 0$. Hence $\nabla^{g}\alpha$ and $\nabla^{g}\beta$ are RadTM-valued since $p$ and $q$ are non zero.

Assume $D_{\alpha}$ is integrable and consider $X,~V \in \Gamma(D_{\alpha})$ and  $Y\in \Gamma(D_{\beta})$. We have

\begin{eqnarray*}
\langle\nabla S(V,Y), X\rangle &=& \langle\nabla_{V}(SY)-S(\nabla_{V}Y),X\rangle\cr
&=& -\beta\langle Y,\nabla_{V}X\rangle-\alpha\langle\nabla_{V}Y,X\rangle\cr
&=& -\beta\langle Y,\nabla_{V}X\rangle + \alpha\langle \nabla_{V}X,Y\rangle\cr
&=&(\alpha - \beta)\langle\nabla_{V}X,Y\rangle.
\end{eqnarray*}
But the last term vanishes since $D_{\alpha}$ is autoparallel from (b) in Fact~\ref{integr}. Thus, we obtain for all $X$, $V$ in $D_{\alpha}$, $Y$ in $D_{\beta}$,
\begin{equation}
\label{eq57}
\langle \nabla S(V,Y), X\rangle =0.
\end{equation}
Now, let $U$ in  $D_{\beta}$. since $C(X,Y)=0$ we have 
\begin{equation}
\label{eq58}
\nabla S(X,Y)=(\beta I-S)\nabla_{X}Y \in S(TM).
\end{equation}
Hence
\begin{equation}
\label{eq59}
\langle\nabla S(X,Y),U\rangle = \langle(\beta I-S)\nabla_{X}Y ,U\rangle = \langle\nabla_{X}Y,(\beta I-S)U\rangle = 0.
\end{equation} 
Similar computation assuming $D_{\beta}$ integrable leads to
\begin{equation}
\label{eq60}
\langle \nabla S(Y,X),U\rangle =0.
\end{equation} 
and 
\begin{equation}
\label{eq61}
\langle \nabla S(Y,X), V\rangle =0.
\end{equation}
for all $Y$, $U$ in $D_{\beta}$ and $X$, $V$ in $D_{\alpha}$. 
Thus, it follows (\ref{eq57}),(\ref{eq59})-(\ref{eq61}) and (b) in Fact~\ref{integr} that $\nabla S$ vanishes on $S(TM)\times S(TM)$.

The converse is immediate from Fact~\ref{integr}.$\square$

\section{Totally umbilic screen foliation}
\label{totumb}
In general a distribution $D\subset TM$ is called umbilical if there exist a vectorfielf $H \in \chi(M)$ such that 
\begin{equation}
\label{62}
\nabla_{X}Y = p(\nabla_{X}Y) + \langle X,X\rangle\varsigma
\end{equation}
for every local section $X\in \Gamma(D)$, where $p$ denotes the "orthogonal" projection $p: TM \longrightarrow D $. In case $D$ is integrable, then it is called totally umbilical. The vector field $\varsigma$ in the definition is called the mean curvature vector of the distribution $D$. In particular, the screen distribution $S(TM)$ is totally umbilical if on any coordinate neighbourhood $\
\mathcal{U}\subset M$there exists a smooth function $\rho$ such that 
\begin{equation}
\label{63}
C(X,PY)= \rho~g(X,Y).
\end{equation}
Now, we state the following

\begin{prop}
\label{umb}
Let $S$ be a $\anabla$-tensor on $(M,g,S(TM))$ where the screen distribution is totally umbilic. Then all the eigenspace distributions $D_{\alpha}=Ker(\alpha I - S)$are umbilical.
\end{prop}

\noindent
{\bf{Proof.~}}
Note that $TM=D_{0} \oplus \sum_{i=1}^{k}D_{i}$ with $D_{0}=RadTM$. Since $\nabla_{\xi}\xi \in D_{0}$, it is obvious that $D_{0}$ is umbilical. Also, for $X\in \Gamma(D_{i})$, 
\begin{equation}
\label{eq64}
\nabla_{X}X = \na_{X}X + C(X,X)\xi = \na_{X}X + \rho g(X,X)\xi 
\end{equation}
Let $p_{i}: TM \longrightarrow D_{i}$ denote the projection morphism on $D_{i}$, we write
$$\nabla_{X}X = p_{i}(\nabla_{X}X) + h_{i}(X,X).$$
It follows that for $Y\in S(TM)$,
$$\langle\nabla_{X}X,Y\rangle = \langle p_{i}(\nabla_{X}X),Y\rangle +\langle h_{i}(X,X),Y\rangle,$$
that is
\begin{eqnarray*}
\langle h_{i}(X,X),Y\rangle &=&\sum_{\stackrel{j=1}{j\ne i}}^{k}\langle\nabla_{X}X,P_{j}Y\rangle\cr 
&\stackrel{(\ref{eq44})}{=}& \frac{1}{2}\langle X,X\rangle\sum_{\stackrel{j=1}{j\ne i}}^{k}\frac{\langle \nabla^{g}\lambda_{i},P_{j}Y\rangle}{\lambda_{j}-\lambda_{i}}\cr
&=&-\frac{1}{2}\langle X,X\rangle\sum_{\stackrel{j=1}{j\ne i}}^{k}\langle \frac{\nabla^{g}\lambda_{i}}{\lambda_{i}-\lambda_{j}},P_{j}Y\rangle\cr
&\stackrel{(\ref{eq43})}{=}&-\frac{1}{2}\langle X,X\rangle\sum_{\stackrel{j=1}{j\ne i}}^{k}\langle P_{j}\nabla^{g}\ln|\lambda_{i}-\lambda_{j}|,Y\rangle.
\end{eqnarray*}
Hence, the $S(TM)$ component of $h_{i}(X,X)$ is 
\begin{equation}
\label{eq65}
\xi_{i}=-\frac{1}{2}\langle X,X\rangle\sum_{\stackrel{j=1}{j\ne i}}^{k}P_{j}(\nabla^{g}\ln|\lambda_{i}-\lambda_{j}|)
\end{equation}
Then, from (\ref{eq64}) we have 
\begin{equation}
\label{eq66}
h_{i}(X,X)= \langle X,X\rangle(\xi_{i}+ \rho \xi).
\end{equation}
Hence $D_{i}$ is umbilical $(1\le i \le k)$, with $\varsigma_{i} = \xi_{i}+\rho \xi$ as mean curvature vector field.

\section{Almost product foliation}
\label{prodstr}
By integrable almost product structure is meant a sequence $(D_{0}, \cdots,D_{k})$ of distributions for which all the distributions $D_{\alpha_{1}}\oplus D_{\alpha_{2}}\oplus \cdots \oplus D_{\alpha_{p}}$ are integrable for any $0\le \alpha_{1}\le \cdots \le \alpha_{p}$ and  $p\in \{0,1,\cdots,k   \}$. A distribution $D_{i}$  ~$(1\le i \le k )$  is called $D_{0}-$almost autoparallel (resp. $D_{0}-$almost parallel) if for any $X$, $Y$ in $\Gamma(D_{i})$, $\nabla_{X}Y \in D_{0}\oplus D_{i} $ (resp. $\forall X \in TM, \forall Y \in \Gamma(D_{i}),~ \nabla_{X}Y \in D_{0}\oplus D_{i}$).

The following result deals with quasi isotropy of $S$. More precisely, we have 

\begin{theo}
\label{t3}
Let $S$ be a $\anabla$-tensor on $(M,g,S(TM))$ with eigenfunctions $\lambda_{0},\lambda_{1},\cdots,\lambda_{k}$. Assume $\nabla^{g}\lambda_{0}$,~$\nabla^{g}\lambda_{1}$, $\dots$, $\nabla^{g}\lambda_{k}$ are $RadTM = D_{0}$- valued and the $D_{\alpha}=Ker(\lambda_{\alpha}I-S)$ define an integrable almost product structure on $M$. Then, $\nabla S|_{S(TM)\times TM} \in \Gamma(D_{0})$.
\end{theo}

\noindent
{\bf{Proof.~}}
First, note that for $X\in S(TM)$, $\nabla S(X,\xi)= (X\cdot\lambda_{0})\xi \in D_{0}$. Now, for $X \in \Gamma(D_{i})$, we have from (\ref{eq42}) and $\nabla^{g}\lambda_{i} \in \Gamma(_{0})$,
\begin{equation}
\label{eq67}
\nabla S(X,X)= (\lambda_{i}-\lambda_{0})~C(X,X)~\xi~ \in \Gamma(D_{0}).
\end{equation}
Integrability of each $D_{i}$ leads to $\nabla S(X,Y)= \nabla S(Y,X)$ for $X$ and $Y$ in $D_{i}$. Also, the integrability of the almost product structure implies $S(TM)$ is integrable and consequently $C$ is symmetric on $S(TM)\times S(TM)$. So, for $X,~Y \in \Gamma(D_{i})$, we obtain by bipolarization of (\ref{eq67}),
\begin{equation}
\label{eq68}
\nabla S(X,Y)= (\lambda_{i}-\lambda_{0})~C(X,Y)~\xi~ \in \Gamma(D_{0}).
\end{equation}
If $X \in \Gamma(D_{i})$ , $Y \in \Gamma(D_{j})$, $i\ne j$, we have from (\ref{eq44}) and $\nabla^{g}\lambda_{i} \in \Gamma(D_{0})$,
\begin{eqnarray*}
\langle\nabla_{X}X,Y\rangle&=&\frac{1}{2}\frac{\langle\nabla^{g}\lambda_{i},Y\rangle }{\lambda_{j}-\lambda_{i}}\langle X,X\rangle ~= ~ 0,
\end{eqnarray*}
Thus, $\nabla_{X}X \in D_{0}\oplus D_{i}$ if $X \in \Gamma(D_{i})$. It follows that 
\begin{eqnarray*}
\nabla_{X}Y + \nabla_{Y}X  \in D_{0}\oplus D_{i}~\mbox{for}~ X ~ \mbox{and} ~Y~\in \Gamma(D_{i}).
\end{eqnarray*}
Since $D_{i}$ is integrable, $\nabla_{X}Y - \nabla_{Y}X \in D_{i}\subset D_{0}\oplus D_{i}$. Hence, for $X$ and $Y$ in $\Gamma(D_{i})$,
\begin{equation}
\label{eq69}
\nabla_{X}Y \in D_{0}\oplus D_{i},
\end{equation}
and each $D_{i}$ is $D_{0}$-almost autoparallel.

Let $i,~j,~l$ be pairwise different numbers and $X \in \Gamma(D_{i})$, $Y \in \Gamma(D_{j})$ and $Z \in \Gamma(D_{l})$. By Koszul formula and integrability of the almost product structure, it follows that
\begin{eqnarray*}
2\langle \nabla_{X}Y,Z\rangle&=& X\cdot\langle Y,Z\rangle + Y\cdot\langle X,Z\rangle- Z\cdot\langle X,Y\rangle \cr & & + \langle [X,Y],Z\rangle +\langle[Z,X],Y\rangle - \langle [Y,Z],X\rangle = 0.
\end{eqnarray*}
Hence, for $X \in \Gamma(D_{i})$,~ $Y \in \Gamma(D_{j})$,~ $(i\ne j)$,
\begin{equation}
\label{eq70}
\nabla_{X}Y \in D_{0}\oplus D_{i}\oplus D_{j}.
\end{equation}
Also, consider $X,~ Z \in \Gamma(D_{i})$,~ $Y \in \Gamma(D_{j})$,~ $(i\ne j)$, we have

\begin{eqnarray*}
0=\langle Z,Y\rangle&\Rightarrow & 0= \langle \nabla_{X}Z,Y\rangle+\langle Z,\nabla_{X}Y\rangle \stackrel{(\ref{eq69})}{=}\langle Z,\nabla_{X}Y\rangle.
\end{eqnarray*}
Then, using (\ref{eq70}) we derive 
\begin{equation}
\label{eq71}
\nabla_{X}Y ~\in D_{0}\oplus D_{j}, ~~\mbox{for} ~X\in D_{i},~ ~Y\in D_{j},~(i\ne j).
\end{equation}
Consequently, from (\ref{eq69}) and (\ref{eq71}), it follows
\begin{equation}
\label{eq72}
\nabla_{X}Y ~\in D_{0}\oplus D_{i}, ~~\mbox{for} ~X\in S(TM) ~\mbox{and}~ Y\in D_{i}.
\end{equation}
Finally, we have from (\ref{eq46}) and (\ref{eq71}) that for $X\in S(TM)$, $Y \in D_{j}$,
\begin{eqnarray}
\label{eq73}
\nabla S(X,Y)&=& \langle \nabla^{g}\lambda_{j},X\rangle Y + (\lambda_{j}I-S)\nabla_{X}Y\cr
& & \cr
&=&(\lambda_{j}I-S)(\eta(\nabla_{X}Y)\xi + p_{j}(\nabla_{X}Y))\cr
& &\cr
&=&(\lambda_{j}-\lambda_{0})~\eta(\nabla_{X}Y)\xi\cr
& &\cr
&=&(\lambda_{j}-\lambda_{0})~C(X,Y)\xi ~ \in \Gamma(D_{0}),
\end{eqnarray}
which completes the proof.$\square$

\begin{rem}
\label{rem2}
It follows that, under hypothesis of theorem~\ref{t3}, we have by (\ref{eq73}),
\begin{equation}
\label{eq74}
\nabla S(X,Y)= (\lambda_{j}-\lambda_{0})~C(X,Y)\xi 
\end{equation} 
for $X\in S(TM)$ and $Y \in D_{j}$.
\end{rem}

\begin{coro}
\label{isotr}
Let $S$ be an $\anabla$-tensor on $(M,g,S(TM))$ with constant eigenfunctions $(\lambda_{0},\lambda_{1},\cdots,\lambda_{k})\in \mathbb{R}^{k+1}$ and integrable almost product structure given by its eigenspace distributions $D_{\alpha}=Ker(\lambda_{\alpha}I-S)$. Then $S$ is an isotropic $\anabla$-tensor. 
\end{coro}

\noindent
{\bf{Proof.~}}From theorem~\ref{t3}, it suffices to show that $\nabla S(\xi,X) \in D_{0}$ for $X$ in $TM$. But since $D_{0}\oplus D_{i}$ is integrable,   $\nabla_{\xi}X \in D_{0}\oplus D_{i}$ for $X \in D_{i}$. Thus, 
\begin{eqnarray*}
\nabla S(\xi,X)&=&(\xi \cdot \lambda_{i})X + (\lambda_{i}-\lambda_{0})p_{0}(\nabla_{\xi}X)\cr
& & \cr
&=&(\lambda_{i}-\lambda_{0})p_{0}(\nabla_{\xi}X) \in D_{0}.\square
\end{eqnarray*}

\begin{coro}
\label{para}
Let $S$ be an $\anabla$-tensor on $(M,g,S(TM))$ with eigenfunctions $(\lambda_{0},\lambda_{1},\cdots,\lambda_{k})$. Assume that  $\nabla^{g}\lambda_{0}$,$\nabla^{g}\lambda_{1}$, $\cdots$, $\nabla^{g}\lambda_{k}$ are $RadTM = D_{0}$- valued and the $D_{\alpha}=Ker(\lambda_{\alpha}I-S)$ define an integrable almost product structure on $M$. If leaves of the screen distribution are totally geodesic in $M$, then $\nabla S = 0$ on $S(TM)\times S(TM)$.
\end{coro}

\noindent {\bf{Proof.~}}The foliation determined by the screen distribution is totally geodesic if and only if $C=0$. Then our claim follows (\ref{eq74}) in remark~\ref{rem2}.$\square$

\section{$\anabla$-tensors and lightlike warped product}
\label{warped}
Lightlike warped products are introduced in \cite{Dug1}, and used in \cite{Dug2} to study the problem of finding globally null manifolds with constant scalar curvature.

Let $(N,g_{N})$and $(F,g_{F})$be a lightlike and a Riemannian manifold of dimension $n$ and $m$ respectively. Let $\pi:N\times F \longrightarrow N$ and $\varrho:N\times F \longrightarrow N$ denote the projection maps given by $\pi(x,y) = x$ and $\varrho (x,y)= y$ for $(x,y)\in N\times F$, respectively, where the projection $\pi$ on $N$ is with respect to a nondegenerate screen distribution $S(TN)$. The product manifold $M=N\times F$, endowed with the degenerate metric defined by 
\begin{equation}
\label{eq75}
g(X,Y)=g_{N}(\pi_{\star}X,\pi_{\star}Y) + f(\pi(x,y))g_{F}(\varrho_{\star}X,\varrho_{\star}Y)
\end{equation}
for all $X$, $Y$ tangent to $M$, where $\star$ is the symbol of the tangent map and $f:N\longrightarrow \mathbb{R}^{\star}_{+}$ is some positive smooth function on $N$. Such a product is denoted $M=(N\times_{f}F,g)$.

\begin{rem}
\label{rem3}
In \cite{Dug1}, this warped product is called of class $A$. The class $B$ one is concerned with  two lightlike factors.
\end{rem}

The following result shows that there is an interplay between existence of $\anabla$-tensors  of certain type and  lightlike warped product structure. In some sense it represents a more general converse to example (a) in section~\ref{examp}.

\begin{theo}
\label{warp}
Let $(M,g)$ be a Killing horizon with a complete simply connected Riemannian global hypersurface,say $H$, $S$ a $\anabla$-tensor on $(M,g,S(TM,\varphi_{t},H))$ with $k+1$ eigenfunctions $\lambda_{0}$~,~$\lambda_{1}$,~$\dots$ ,~$\lambda_{k}$ and eigendistributions $D_{0}=Ker(S-\lambda_{0}I)=RadTM$, $D_{i}=Ker(S-\lambda_{i}I),~~i=1,\dots,k$. If
\begin{enumerate}
\item[(a)]$\lambda_{1}=\mu = \mbox{constant}$,
\item[(b)]The almost product structure $(D_{0},D_{1},\dots,D_{k})$ is integrable,
\item[(c)]$\nabla\lambda_{i} \in D_{1},~~i=1,\dots,k$
\end{enumerate}
then,
\[M = \mathbb{L}\times M_{1}\times_{f_{2}}M_{2}\times \cdots,\times_{f_{k}}M_{k}  \]
where $\mathbb{L}$ (a global null curve )is a one-dimensional integral manifold of the global null vector on $M$ and $M_{i}~~(1\le i \le k)$ are leaves of $D_{i}$ and $f_{i}^{2}=|\lambda_{i}-\mu|$, $(2\le i \le k)$ are smooth positive functions on $M_{1}$.
\end{theo}

\noindent
{\bf{Proof.~}}First, note that thr radical distribution of $M$ is spaned by a global null killing vector field and the screen distribution $S(TM)$ is integrable. It follows that $M = \mathbb{L} \times M'$ is a global product manifold where  $\mathbb{L}$ is a one-dimensional integral manifold of a global null vector field on $M$, and $M'$ a leaf of $S(TM)$. Let $g'$ denote the Riemannian metric induced on $M'$ and $\pi:\mathbb{L}\times M'\longrightarrow M'$ the natural projection map. Then the lightlike hypersurface $(M,g)$  is isommetric to $(\mathbb{L}\times M', g=\pi^{\star}g')$. Also, by (b) in remark~\ref{rem1bis},$S$ induces by restriction an $\mathcal{A}$-tensor $S'$ on $M'$ with respect to the Levi-Civita connection $\stackrel{\star}{\nabla}$ it inherits from $M$ and $\lambda_{1},\dots,\lambda_{k}$ are eigenfunctions of $S'$, with eigendistributions $D_{i}$, $i=1,\dots,k$. Since $\nabla \lambda_{i} \in D_{1}$,~$i=1,\dots,k$, the functions $\lambda_{i}$ depend uniquely on the parameters on the leaf $M_{1}$ of the integrable almost product structure $(D_{0},D_{1},\dots,D_{k})$. In particular, 
$$\oplus_{j>1}D_{j}\subset Kerd\lambda_{i},~~(1\le i \le k).$$
The final result follows \cite{Jel2}. Indeed, since in addition to above facts, $H$ (and then $M'$) is complete simply connected Riemannian hypersurface of $M$, we have

$$(M',g') = M_{1}\times_{f_{2}}M_{2}\times \dots \times_{f_{k}}M_{k} $$
where $TM_{i}= D_{i}$~ and~ $f_{i}=\sqrt{|\lambda_{i}-\mu|}$, ~$2\le i \le k$. Then 
\[ M= \mathbb{L}\times M_{1}\times_{f_{2}}M_{2}\times \dots \times_{f_{k}}M_{k} \]
is a multiply warped product manifold where $f_{2},\dots,f_{k}$are smooth positive functions on the factor $M_{1}$ of the lightlike product manifold $\mathbb{L}\times M_{1}$.$\square$

\vspace{1.5cm}
{\bf{Acknowledgments.}}
The first named author (C. Atindogbe) thanks the Agence Universitaire de la Francophonie (AUF) for support with a one year research grant, along with the Institut Elie Cartan (IECN, UHP-Nancy~I) for  research facilities during the completion of this work.

\end{document}